\documentclass[12pt]{article}

\usepackage{amssymb}

\newtheorem{theorem}{$~~~~$ Theorem}[section]

\newtheorem{example}[theorem]{$~~~~$ Example}

\newtheorem{corollary}[theorem]{$~~~~$ Corollary}
\newtheorem{lemma}[theorem]{$~~~~$ Lemma}
\newtheorem{remark}[theorem]{$~~~~$Remark}
\newtheorem{definition}[theorem]{$~~~~$Definition}

\newtheorem{remm}[theorem]{$~~~~$Remark}

\newtheorem{deff}[theorem]{$~~~~$Definition}

\def\ulxyz{\ulcorner}
\def\urxyz{\urcorner}

\def\aaa{\beta}

\def\nor1{Normed$\{~2^{ \zzz \theta  \, )} ~$,$~\sqrt{~2^{ \zzz \theta  \, )}}~\}$}

\def\xor2{Normed$\{ ~\sqrt{~2^{ \zzz \theta  \, )}}~,~2~ \} $}

\def\pag2{Page 2}

\def\zzz{~ \sharp ( ~ }

\def\f55{ \normalsize  \baselineskip = 1.8 \normalbaselineskip }

\def\f55{  \baselineskip = 1.1 \normalbaselineskip } 
\def\g55{  \baselineskip = 1.0 \normalbaselineskip } 
\def\s55{ \baselineskip = 1.0 \normalbaselineskip } 

\def\f55{  \baselineskip = 0.7 \normalbaselineskip }

\newcommand{\thx}[1]{Theorem \ref{#1}}

\newcommand{\el}[1]{Line (\ref{#1})}

\newcommand{\eq}[1]{(\ref{#1})}

\begin{document}


%



\title{How the Law of Excluded Middle
Pertains to
the Second Incompleteness Theorem
and its Boundary-Case Exceptions}


\def\beq{\begin{equation}}
\def\enq{\end{equation}}

\def\bel{\begin{lemma}}
\def\enl{\end{lemma}}

\def\bec{\begin{corollary}}
\def\enc{\end{corollary}}

\def\bed{\begin{description}}
\def\ennd{\end{description}}
\def\bee{\begin{enumerate}}
\def\ene{\end{enumerate}}

\def\bxbxd{\begin{definition}}
\def\bxbxdd{\begin{definition}}
\def\eedd{\end{definition}}
\def\bxbxdr{\begin{definition} \rm}
\def\bel{\begin{lemma}}
\def\enl{\end{lemma}}
\def\ent{\end{theorem}}

\author{  Dan E.Willard} 




\date{State University of New York at Albany}

\maketitle

 \setcounter{page}{0}

 \thispagestyle{empty}



\normalsize

\baselineskip = 1.3\normalbaselineskip

\normalsize

\baselineskip = 1.0 \normalbaselineskip 
\def\bbint{\large \baselineskip = 1.6 \normalbaselineskip } 
\def\bbint{\large \baselineskip = 1.6 \normalbaselineskip }
\def\bbint{\normalsize \baselineskip = 1.3 \normalbaselineskip }


\def\bbint{\normalsize \baselineskip = 1.27 \normalbaselineskip }

\def\bbint{\large \baselineskip = 2.0 \normalbaselineskip }

\def\bbint{\normalsize \baselineskip = 1.25 \normalbaselineskip }
\def\bbina{\normalsize \baselineskip = 1.24 \normalbaselineskip }

\def\bbint{\large \baselineskip = 2.0 \normalbaselineskip }

\def\bbing{ }
\def\bbins{ }
\def\bbinm{ }

\def\bbint{\normalsize \baselineskip = 1.95 \normalbaselineskip }

\def\bbing{ }
\def\bbins{ }
\def\bbinm{ }

\def\bbint{\large \baselineskip = 2.3 \normalbaselineskip } 
\def\bbing{ }
\def\bbins{ }
\def\bbinm{ }

\def\bbint{\normalsize \baselineskip = 1.7 \normalbaselineskip } 

\def\bbint{\large \baselineskip = 2.3 \normalbaselineskip } 
\def\bbinm{ \baselineskip = 1.18 \normalbaselineskip }

\def\bbint{\large \baselineskip = 2.0 \normalbaselineskip } 
\def\bbing{ }
\def\bbins{ }
\def\bbinm{ }
\def\bbinr{ }

\def\bbint{\normalsize \baselineskip = 1.25 \normalbaselineskip }
\def\bbina{\normalsize \baselineskip = 1.24 \normalbaselineskip }
\def\bbinr{ \baselineskip = 1.3 \normalbaselineskip }
\def\bbing{ \baselineskip = 1.28 \normalbaselineskip }
\def\bbins{ \baselineskip = 1.21 \normalbaselineskip }
\def\bbinm{  }

\def\ftl{ \baselineskip = 1.5 \normalbaselineskip }

\bbint

\noindent

\small

\baselineskip = 1.14 \normalbaselineskip

\baselineskip = 1.2 \normalbaselineskip 



\large
\normalsize

 \baselineskip = 1.2 \normalbaselineskip 


\begin{abstract}
\large
\baselineskip = 1.5 \normalbaselineskip  

Our earlier publications showed semantic tableau admits partial exceptions to
the Second Incompleteness Theorem where a formalism recognizes its self
consistency and views multiplication as a 3-way relation (rather than as a
total function).  We now show these boundary-case evasions will
collapse if the Law
of the Excluded Middle is treated by tableau as a schema of logical axioms
(instead of as derived theorems).

\end{abstract}  

  \baselineskip = 1.3 \normalbaselineskip

\bigskip


\bigskip

\bigskip

\bigskip

\bigskip


\bigskip

\bigskip
 
 \bigskip

 
 \bigskip
 \bigskip
 \bigskip

 \small
 {\bf Comment:$~$ } 
 Short conference announcements of these results 
at ASL-2020's Virtual N. American Meeting and at
LFCS-2020.

\def\ww22{\normalsize \baselineskip = 1.21\normalbaselineskip \parskip 4 pt}
\def\bb22{\normalsize \baselineskip = 1.19\normalbaselineskip \parskip 4 pt}
\def\zz22z{\normalsize \baselineskip = 1.19 \normalbaselineskip \parskip 3 pt}
\def\xx22{\normalsize \baselineskip = 1.17\normalbaselineskip \parskip 4 pt}
\def\vx22s{\normalsize \baselineskip = 1.16 \normalbaselineskip \parskip 3 pt} 
\def\vv22{\normalsize \baselineskip = 1.17 \normalbaselineskip \parskip 3 pt} 
\def\aa22{\normalsize \baselineskip = 1.15 \normalbaselineskip \parskip 3 pt} 
\def\g55{  \baselineskip = 1.0 \normalbaselineskip } 
\def\s55{ \baselineskip = 1.0 \normalbaselineskip } 
\def\sm55{ \baselineskip = 0.9 \normalbaselineskip }

\vspace*{- 1.0 em}

\def\waw11{\normalsize \baselineskip = 1.72\normalbaselineskip}
\def\waw11{\normalsize \baselineskip = 1.12\normalbaselineskip}
\def\waw11{\normalsize \baselineskip = 1.85\normalbaselineskip}

\def\waw11{\normalsize \baselineskip = 1.45\normalbaselineskip}

\def\waw11{\normalsize \baselineskip = 1.7\normalbaselineskip}

\def\waw11{\normalsize \baselineskip = 1.12\normalbaselineskip}

\def\g55{  \baselineskip = 1.50 \normalbaselineskip } 
\def\s55{ \baselineskip = 1.50 \normalbaselineskip } 
\def\sm55{ \baselineskip = 1.5 \normalbaselineskip }

\def\g55{  \baselineskip = 1.50 \normalbaselineskip } 
\def\s55{ \baselineskip = 1.50 \normalbaselineskip } 
\def\sm55{ \baselineskip = 0.9 \normalbaselineskip }

\def\aa22{\normalsize  \waw11 \parskip 6 pt} 
\def\bb22{\normalsize  \waw11 \parskip 5 pt}
\def\ww22{\normalsize \waw11 \parskip 4 pt}
\def\vv22{\normalsize  \waw11 \parskip 3 pt} 
\def\tt22{\normalsize  \waw11 \parskip 2 pt} 

\def\g55{  \baselineskip = 1.0 \normalbaselineskip } 
\def\b55{  \baselineskip = 1.0 \normalbaselineskip } 
\def\s55{ \baselineskip = 1.0 \normalbaselineskip } 
\def\sm55{ \baselineskip = 0.9 \normalbaselineskip }

\def\mal{ \bf  }
\def\nal{\mathcal}

\def\cvrew{ \baselineskip = 1.6 \normalbaselineskip \parskip 3pt }

\def\ttt2c{ }
\def\tttc{ }

\def\tttc{\tiny \baselineskip = 0.8 \normalbaselineskip  \parskip 0pt }
\def\ttt2c{\tiny \baselineskip = 0.7 \normalbaselineskip  \parskip 0pt }
\def\tttc{ \baselineskip = 2.1 \normalbaselineskip  \parskip 5pt }
\def\ttt2c{ \baselineskip = 2.1 \normalbaselineskip  \parskip 5pt }

\def\tttc{ \baselineskip = 1.15 \normalbaselineskip  \parskip 5pt }
\def\ttt2c{ \baselineskip = 1.15 \normalbaselineskip  \parskip 5pt }

\def\tttc{ \baselineskip = 1.12 \normalbaselineskip  \parskip 4pt }
\def\ttt2c{ \baselineskip = 1.12 \normalbaselineskip  \parskip 4pt }

\def\tttc{ \baselineskip = 1.14 \normalbaselineskip  \parskip 3pt }
\def\ttt2c{ \baselineskip = 1.10 \normalbaselineskip  \parskip 2pt }
\def\ttt2c{ \baselineskip = 0.98 \normalbaselineskip  \parskip 0pt }

\def\cvt{ \baselineskip = 0.98 \normalbaselineskip }
\def\cv9{ \baselineskip = 0.99 \normalbaselineskip }
\def\cvs{ \baselineskip = 1.0 \normalbaselineskip }
\def\cvl{ \baselineskip = 1.0 \normalbaselineskip }
\def\cvh{ \baselineskip = 1.03 \normalbaselineskip }
\def\cvg{ \baselineskip = 1.00 \normalbaselineskip }

\def\cvt{ \baselineskip = 1.6 \normalbaselineskip }
\def\cv9{ \baselineskip = 1.6 \normalbaselineskip }
\def\cvs{ \baselineskip = 1.6 \normalbaselineskip }
\def\cvl{ \baselineskip = 1.6 \normalbaselineskip }
\def\cvh{ \baselineskip = 1.6 \normalbaselineskip }
\def\cvg{ \baselineskip = 1.6 \normalbaselineskip }
\def\cvb{ \baselineskip = 1.6 \normalbaselineskip }
\def\cvnew{ \baselineskip = 1.6 \normalbaselineskip }
\def\cvmew{ \baselineskip = 1.6 \normalbaselineskip }
\def\cvwew{ \baselineskip = 1.6 \normalbaselineskip \parskip 5pt }
\def\cvrew{ \baselineskip = 1.6 \normalbaselineskip \parskip 3pt }

\def\cvt{ \baselineskip = 1.22 \normalbaselineskip }
\def\cv9{ \baselineskip = 1.22 \normalbaselineskip }
\def\cvs{ \baselineskip = 1.22 \normalbaselineskip }
\def\cvl{ \baselineskip = 1.22 \normalbaselineskip }
\def\cvh{ \baselineskip = 1.22 \normalbaselineskip }
\def\cvg{ \baselineskip = 1.22 \normalbaselineskip }
\def\cvb{ \baselineskip = 1.22 \normalbaselineskip }
\def\cvnew{ \baselineskip = 1.4 \normalbaselineskip }
\def\cvmew{ \baselineskip = 1.35 \normalbaselineskip }
\def\cvwew{ \baselineskip = 1.4 \normalbaselineskip \parskip 5pt }
\def\cvrew{ \baselineskip = 1.22 \normalbaselineskip \parskip 3pt }

\def\cvt{ }
\def\cv9{ }
\def\cvs{ }
\def\cvl{ }
\def\cvh{ }
\def\cvg{ }
\def\cvb{ }
\def\cvnew{ } 
\def\cvmew{ }
\def\cvwew{ }
\def\cvrew{ }

\def\fend{ 

\medskip -------------------------------------------------------}

\def\g55{  \baselineskip = 1.0 \normalbaselineskip } 
\def\s55{ \baselineskip = 1.0 \normalbaselineskip } 
\def\sm55{ \baselineskip = 1.0 \normalbaselineskip } 
\def\h55{  \baselineskip = 1.08 \normalbaselineskip } 
\def\b55{  \baselineskip = 1.1 \normalbaselineskip } 

\normalsize

\baselineskip = 1.85 \normalbaselineskip





\vspace*{- 1.0 em}



\baselineskip = 1.04 \normalbaselineskip 

\baselineskip = 0.96 \normalbaselineskip 

%
\baselineskip = 2.16 \normalbaselineskip 
\baselineskip = 2.3 \normalbaselineskip 

\baselineskip = 0.95 \normalbaselineskip 
\baselineskip = 0.88 \normalbaselineskip 
\parskip 0pt

\noindent

%
%
%
%


\def\gvs{ \normalsize \baselineskip = 1.4 \normalbaselineskip  \parskip    5pt}
\def\gvs{ \normalsize \baselineskip = 1.44 \normalbaselineskip  \parskip    5pt}
\def\gvs{ \large \baselineskip = 1.44 \normalbaselineskip  \parskip    5pt}
\def\gvs{ \normalsize \baselineskip = 1.44 \normalbaselineskip  \parskip    5pt}\def\gvs{ \normalsize \baselineskip = 1.74 \normalbaselineskip  \parskip    5pt}
\def\gvs{ \normalsize \baselineskip = 1.44 \normalbaselineskip  \parskip 5pt}

\def\gvs{   \baselineskip = 1.74 \normalbaselineskip  \parskip    5pt}

\def\gvs{ \normalsize \baselineskip = 1.44 \normalbaselineskip  \parskip 5pt}
\def\gvs{ \large \baselineskip = 2.0 \normalbaselineskip  \parskip 5pt}
\def\gvs{ \Large \baselineskip = 2.0 \normalbaselineskip  \parskip 5pt}
\def\gvs{ \normalsize \baselineskip = 2.44 \normalbaselineskip  \parskip 5pt}
\def\gvs{ \normalsize \baselineskip = 2.04 \normalbaselineskip  \parskip 5pt}
\def\gvs{ \normalsize \baselineskip = 2.64 \normalbaselineskip  \parskip 5pt}
\def\gvs{ \Large \baselineskip = 1.6 \normalbaselineskip  \parskip 5pt}

\gvs

\footnotesize

\def\gvs{ }

\normalsize \baselineskip = 0.98 \normalbaselineskip
\normalsize \baselineskip = 1.0 \normalbaselineskip
\normalsize \baselineskip = 1.01 \normalbaselineskip

\def\gvs{ \normalsize \baselineskip = 1.25 \normalbaselineskip  \parskip 4pt}

\def\gvs{ \Large \baselineskip = 1.6  \normalbaselineskip  \parskip 6pt}
\def\gvs{ \normalsize \baselineskip = 1.6  \normalbaselineskip  \parskip 6pt}
\def\gvs{ \large \baselineskip = 1.6  \normalbaselineskip  \parskip 6pt}

\def\gvs{ \normalsize \baselineskip = 1.227 \normalbaselineskip  \parskip 3pt}
\def\gvs{ \large \baselineskip = 1.8  \normalbaselineskip  \parskip 6pt}

\def\gvs{ \normalsize \baselineskip = 1.5 \normalbaselineskip  \parskip 3pt}

\def\gvs{ \large \baselineskip = 2.1  \normalbaselineskip  \parskip 6pt}

\def\gvs{ \normalsize \baselineskip = 2.1  \normalbaselineskip  \parskip 6pt}


 \def\gvs{ \normalsize \baselineskip = 1.227 \normalbaselineskip  \parskip 3pt}

 \def\gvs{ \large  \baselineskip = 1.6 \normalbaselineskip  \parskip 5pt}

\def\gvs{ \Large  \baselineskip = 1.8 \normalbaselineskip  \parskip 5pt}
\def\gvs{ \LARGE  \baselineskip = 1.8 \normalbaselineskip  \parskip 5pt}
\def\gvs{ \normalsize  \baselineskip = 2.0 \normalbaselineskip  \parskip 5pt}

\def\gvs{ \Large  \baselineskip = 2.0 \normalbaselineskip  \parskip 5pt}

\def\gvs{ \large  \baselineskip = 2.2 \normalbaselineskip  \parskip 5pt}

\def\gvs{ \normalsize \baselineskip = 2.4  \normalbaselineskip  \parskip 6pt}

\def\gvs{ \normalsize \baselineskip = 2.6  \normalbaselineskip  \parskip 6pt}
\def\gvs{ \normalsize \baselineskip = 2.2  \normalbaselineskip  \parskip 6pt}
\def\gvs{ \normalsize \baselineskip = 1.8  \normalbaselineskip  \parskip 5pt}
\def\tttc{ }

\def\gv2{ \normalsize \baselineskip = 1.30  \normalbaselineskip  \parskip 3pt}

\def\gvs{ }

\def\gvs{ \normalsize \baselineskip = 2.1 \normalbaselineskip  \parskip 7pt}
\def\gvs{ \normalsize \baselineskip = 1.8 \normalbaselineskip  \parskip    7pt}

 \def\gvs{ \large \baselineskip = 1.7  \normalbaselineskip  \parskip 9pt}
\def\gvs{ \normalsize \baselineskip = 2.0  \normalbaselineskip  \parskip 9pt}

\def\gv2{ \normalsize \baselineskip = 1.30  \normalbaselineskip  \parskip 3pt}

\def\gvs{ \large \baselineskip = 1.7  \normalbaselineskip  \parskip 5pt}

\def\gvs{ \normalsize \baselineskip = 2.0  \normalbaselineskip  \parskip 8pt}
\def\gvs{ \large \baselineskip = 2.0  \normalbaselineskip  \parskip 8pt}


\def\gvx{ \large \baselineskip = 1.6 \normalbaselineskip  \parskip    3pt}

\def\gvs{ \Large \baselineskip = 1.8 \normalbaselineskip  \parskip   5pt}

\def\gvs{ \LARGE \baselineskip = 2.0 \normalbaselineskip  \parskip   5pt}

\def\gvs{ \normalsize \baselineskip = 2.0 \normalbaselineskip  \parskip   5pt}
\def\gvs{ \normalsize \baselineskip = 1.4 \normalbaselineskip  \parskip   5pt}

 \def\gvs{ \large \baselineskip = 2.0 \normalbaselineskip  \parskip   5pt}
 \def\gvs{ \Large \baselineskip = 2.0 \normalbaselineskip  \parskip   5pt}
\def\gvs{ \normalsize \baselineskip = 2.0 \normalbaselineskip  \parskip   5pt}
\def\gvs{ \normalsize \baselineskip = 2.5 \normalbaselineskip  \parskip   5pt}

\def\gvs{ \normalsize \baselineskip = 1.4 \normalbaselineskip  \parskip   5pt}
\def\gvs{ \large      \baselineskip = 1.7 \normalbaselineskip  \parskip   5pt}
\def\gvs{ \normalsize \baselineskip = 1.7 \normalbaselineskip  \parskip   5pt}

\def\gvs{ \large \baselineskip = 1.6 \normalbaselineskip}

\def\gvs{ \normalsize \baselineskip = 1.8 \normalbaselineskip}

\def\gvs{ \Large \baselineskip = 1.8 \normalbaselineskip}

\def\gvs{ \normalsize \baselineskip = 2.2 \normalbaselineskip}
\def\gvs{ \normalsize \baselineskip = 2.7 \normalbaselineskip}

\def\gvs{ \normalsize \baselineskip = 2.9 \normalbaselineskip}
\def\gvs{ \normalsize \baselineskip = 3.0 \normalbaselineskip}

\def\gvs{ \normalsize \baselineskip = 2.0 \normalbaselineskip}
\def\gvs{ \normalsize \baselineskip = 1.7 \normalbaselineskip}
\def\gvs{ \normalsize \baselineskip = 1.4 \normalbaselineskip}

\def\gvs{ \normalsize \baselineskip = 1.0 \normalbaselineskip}
\def\gvs{ \large \baselineskip = 1.7 \normalbaselineskip}

\def\gvs{ \normalsize \baselineskip = 1.0 \normalbaselineskip}

\def\gvs{ \normalsize \baselineskip = 1.89 \normalbaselineskip}

\def\gvs{ \large \baselineskip = 1.7 \normalbaselineskip}

\def\gvs{ \normalsize \baselineskip = 2.0 \normalbaselineskip}
\def\gvs{ \large \baselineskip = 2.0 \normalbaselineskip}

\def\gvs{ \normalsize \baselineskip = 1.4 \normalbaselineskip}

\def\gvs{ \large \baselineskip = 2.0 \normalbaselineskip}

\def\gvs{ \normalsize \baselineskip = 2.7 \normalbaselineskip}

\def\gvs{ \normalsize \baselineskip = 1.37 \normalbaselineskip}

\def\gvs{ \normalsize \baselineskip = 1.4 \normalbaselineskip}





\newpage

\section{Introduction }

\label{nnn1}
\gvs
\tttc

This article is intended to explore
the {\it ``hidden significance''} and unexplored
implications of
G\"{o}del's Second Incompleteness Theorem
and
its
 various
generalizations.
In particular,
the existence of a deep
 chasm separating 
the  goals of 
Hilbert's consistency program from the implications of
the
Second Incompleteness Theorem was
 evident,
  immediately,
after 
G\"{o}del
published
 \cite{Go31}'s seminal
announcement.
We 
exhibited
in  \cite{ww93,ww1,ww2,ww5,wwapal,ww6,ww7,ww9}
 a large number of
articles
 about generalizations and boundary case
exceptions
 to the Second Incompleteness Theorem, starting 
with
our  1993 article \cite{ww93}.
These papers, which 
included six papers
published
 in the JSL and APAL,
showed 
every
 extension $~\alpha~$ of Peano Arithmetic can be
mapped onto an axiom system $~\alpha^*$ 
that can
recognize its own consistency and
 prove
analogs of all  $~\alpha$'s $\Pi_1$ 
theorems
(in a slightly different
language, called $~L^*$).

The term 
 ``Self-Justifying'' arithmetic
was employed in our  articles.
 \cite{ww1,ww5,wwapal,ww6,ww9}.
These papers were able
to verify their own consistency by containing a
built-in self-referencing
 axiom 
that
 declared {\it  ``I am consistent''}
(as will be  explained later).
In particular, 
our axiom systems 
 $~\alpha^*~$ 
used the
Fixed-Point Theorem 
to assure
 $~\alpha^*\,$'s 
self-referencing
 analogs of
the pronoun ``I'' would enable 
it
 to refer to 
itself
 in the context of its
{\it  ``I am consistent''} axiomatic declaration.


It turns out that such a self-referencing mechanism will produce
unacceptable G\"{o}del-style
diagonalizing
 contradictions, when
either
 $\alpha^*$ or
its 
particular
deployed
 definition of 
consistency are
too
 strong.
This is because our
           methodologies
{\it only}
    become contradiction-free {\it  when }
 $\alpha^*$ 
uses 
sufficiently
       weak  underlying
         structures.

These
weak structures obviously have
significant
 disadvantages.
Their virtue is that
their
formalisms
  $\alpha^*$ can 
be arranged to
prove more $\Pi_1$ like
theorems than
 Peano Arithmetic,
while offering
{\it 
some type}
 of
{\it partial}
 knowledge about their 
 own consistency.
We will call such formalisms   {\bf ``Declarative Exceptions''}
to the Second Incompleteness Theorem.

An
alternative type of
exception to the Second
Incompleteness Theorem,
which we shall
 call  an {\bf ``Infinite-Ranged
 Exception''},
 was recently developed
by 
Sergei
Artemov  \cite{Ar19}
(It is  related to the
 works of
Beklemishev 
\cite{Be5}
and Artemov-Beklemishev \cite{AB5}.) 
Artemov 
observed Peano Arithmetic 
can 
verify its own consistency, from
a special infinite-ranging perspective.
This means 
PA will generate an infinite set of theorems
$T_1,~T_2,~T_3~...~$ 
where each  $~T_i~$ shows some subset
$~S_i~$
 of PA is unable to prove $0=1$ and
where PA equals
 the formal
union of
these
 special
 selected
 $~S_i~$
satisfying the inclusion property of
 $~S_1 \,  \subset  \, S_2 \,  \subset  \, S_3 \, \subset~ ...~$.

This perspective,
which is
certainly very
useful,
 is
also
 not a panacea. Thus,  the
abstract in
  \cite{Ar19}
cautiously
used
 the adjective of ``somewhat'' to describe how
it  sought to
partially achieve the goals sought by Hilbert's
Consistency Program (with an infinite collection of theorems
$T_1,~T_2,~T_3~...~$ 
 replacing Hilbert's 
intended
goal of
finding 
one unifying formal consistency theorem).

Our
 ``Declarative'''  exceptions
 to the Second Incompleteness
Theorem
and Artemov's  ``Infinite Ranging'' exceptions
are
two quite different
rigorous results, which are
nicely
compatible
 with each other.
This is 
 because each
 acknowledged
that the Second Incompleteness Theorem is a strong
result,  that 
{\it will
admit no full-scale
exceptions.}
Also,
these results are
of interest
because
G\"{o}del 
openly 
conjectured
that
Hilbert's Consistency Program
would
ultimately, reach
{\it some 
 levels of
 partial
 success}
(see next section).
We will explain, 
herein,
 how G\"{o}del's conjecture can be
{\it partially justified,} 
due to an
 unusual consequence of the Law of the Excluded Middle.



More specifically,
we shall focus on
the
semantic tableau deductive mechanisms
of Fitting and Smullyan 
\cite{Fi96,Sm95}
and 
their
special properties from the perspective of
our JSL-2005 article \cite{ww5}.
Each instance of the Law of the Excluded Middle
has been  treated by
most
 tableau mechanisms
 as a provable
theorem,  rather than  as a built-in logical axiom.
This may, at first, appear to be an
 insignificant
distraction
because most
deductive methodologies do not have their consistency
reversed when a theorem 
 is promoted into becoming a logical
axiom.

Our self-justifying axiom systems are
{\it  different,}
however, because their built-in self-referencing {\it ``I am consistent"} axioms
have  their meanings 
change,
 fundamentally,
 when 
their self-referencing
concept of ``I''
involves
promoting
a schema of theorems verifying the Law of Excluded Middle
{\it into
formal explicitly declared
logical axioms.}

This 
effect
is counterintuitive
because similar  distinctions  exist almost nowhere
 else  in Logic.
Thus some confusion, that has surrounded our prior
work, can be clarified when one realizes that
{\it an  interaction}
between
 the self-referencing
concept of ``I'' with 
the  Law of Excluded Middle
causes
the Second Incompleteness
Theorem
 to become activated
{\it precisely
when} the  Law of Excluded Middle
{\it is promoted} 
into becoming a 
schema of
logical axioms.

The intuitive reason for this
unusual
 effect is that
the
transforming
of derived theorems  {\it into}
 logical axioms 
can shorten proofs under the Fitting-Smullyan
semantic tableau
technology.
In the
 particular
 context where 
\textsection \ref{s3}'s  formalism
uses
 self-referencing {\it ``I am consistent"} axioms
and   views
 multiplication as a 3-way relationship,
these conditions will be sufficient for
enacting the full power of the
Second Incompleteness Theorem. 


The next chapter will  explain
how 
these issues are  related to 
questions raised by 
 G\"{o}del and
Hilbert
about
 feasible
 boundary-case exceptions to
the   Second Incompleteness
Effect.


\section{Revisiting
Some Intuitions of
G\"{o}del
and Hilbert}

\label{nnn2}
\gvs
\tttc

Interestingly,
neither  G\"{o}del
 (unequivocally) 
nor Hilbert (after learning about 
 G\"{o}del's work)
would dismiss the possibility of a
compromise
solution, whereby a
{\it  fragment} of the goals of Hilbert's
Consistency Program would remain intact.
Thus, Hilbert never withdrew
 \cite{Hil26}'s
 statement $*$ 
for justifying
his program:
\begin{quote}
\small
\baselineskip = 1.0 \normalbaselineskip
\ttt2c 
$*~$
{\it ``
Let us admit that the situation in which we presently
find ourselves with respect to paradoxes is in the long
run intolerable. Just think: in mathematics, this paragon of
reliability and truth, the very notions and inferences,
as everyone learns, teaches, and uses them, lead to absurdities.
And 
where 
else 
would 
reliability and truth be found 
if  even mathematical thinking fails?''}
\end{quote}
%


G\"{o}del was, also, 
 cautious (especially during the early 1930's) not
to speculate
 whether all facets of Hilbert's Consistency program
would come to a termination.
He thus inserted the following
cautious
caveat into
his
famous 1931 
paper  \cite{Go31}:
\begin{quote}
\small
\ttt2c 
$~**~~$ 
{\it ``It must be expressly noted that
Theorem XI''}
(e.g. the Second Incompleteness Theorem) 
{\it ``represents no contradiction of the formalistic
standpoint of Hilbert. For this standpoint
presupposes only the existence of a consistency
proof by finite means, and {\it there might
conceivably be finite proofs} which cannot
be stated in P or in ... ''}
\end{quote}
Several 
biographies
of 
 G\"{o}del
\cite{Da97,Go5,Yo5}
have noted that
 G\"{o}del's 
 intention (prior to 1930)
was
to
establish
Hilbert's proposed objectives, before
he formalized
his famous result 
that led
in
an
opposite direction.
Moreover,
 Yourgrau's
biography
 \cite{Yo5}
 of G\"{o}del 
records
 how
von Neumann 
found it necessary during the early 1930's to
{\it ``argue 
against G\"{o}del 
himself''}
 about the definitive 
 termination of Hilbert's
consistency program,  
which
{\it ``for several years''} after \cite{Go31}'s publication,
G\"{o}del 
{\it ``was cautious not to prejudge''}.

It is known that
 G\"{o}del 
hinted
the Second Incompleteness
Theorem was more significant 
in
 a 1933 Vienna
lecture  \cite{Go33}.
Yet,  G\"{o}del 
(who  published only
about
85 pages
 during his career) was 
frequently 
 ambivalent 
about this point.
Thus,
a YouTube talk by 
 Gerald Sacks  \cite{YouSa14}
recalled G\"{o}del
telling 
Sacks
some type of 
revival of Hilbert's Consistency
Program
was
likely
 (see
 footnote 
\footnote{\label{f2}
Some quotes from Sacks's
YouTube talk
\cite{YouSa14} are that G\"{o}del
 {\it ``did not  think''}
the objectives of Hilbert's Consistency Program 
{\it ``were erased''} 
by
the Incompleteness Theorem, and
G\"{o}del believed (according to Sacks) 
 it left
  Hilbert's program 
{\it ``very much alive and
even more interesting than it initially was''}. 
} for more details).
Moreover, Anil  Nerode
 has  told us
\cite{Ne20}
he recalled
Stanley Tennenbaum
 having similar conversations
with G\"{o}del,
where G\"{o}del
again stated his 
suspicion
 that  Hilbert's
Consistency Program would be partially revived.
Many scholars have been caught
 by surprise by
G\"{o}del's private
hesitation
about the broader implications of
the Second  Incompleteness Effect.
This is because  G\"{o}del 
only
published roughly 85 pages
during his career, and he
never 
publicly
expanded upon  \cite{Go31}'s 
statement ** .

The research that followed
 G\"{o}del's  seminal  1931 discovery has 
technically
 focused
on studying
mostly
 generalizations of the Second Incompleteness
Theorem
(instead of also examining its
 boundary-case
exceptions). Many of  these generalizations 
of the Second Incompleteness Theorem 
\cite{AZ1,Ar1,Be14,BS76,Bu86,BI95,Fe60,Fr79a,Ha7,Ha11,HP91,KT74,Pa71,PD83,Pa72,Pu85,Pu96,So88,So94,Sv7,Vi5,WP87,ww1,ww2,wwlogos,wwapal}  
are quite subtle.

    The author of this paper is
            especially impressed by a generalization of the
Second Incompleteness Effect, arrived at by the 
combined work of Pudl\'{a}k and  Solovay
together with
added research by
 Nelson and Wilkie-Paris
\cite{Ne86,Pu85,So94,WP87}.
These results, which 
  have been
further
amplified
in \cite{BI95,Fr79a,Ha7,Sv7,ww1},
show
            the Second Incompleteness Theorem
does not require the presence of the Principle of Induction
to apply to
 most 
formalisms that use a Hilbert-Frege 
style
 of 
deduction. 

The
next chapter's Remark \ref{nremm-2.5}
will 
helpfully
 summarize such
 generalizations
 of
the Second Incompleteness Effect.

\section{ Main Notation and Background Literature}

\label{s3}
\label{nnn3}

Let us 
call an
ordered pair $(\alpha,D)$ a
    {\bf Generalized Arithmetic Configuration}
(abbreviated as a {\bf ``GenAC'' })
when  its 
first and second 
components 
are 
defined 
as 
follows:
\bee
\item
The {\bf Axiom Basis} ``$~\alpha~$'' 
for a 
 GenAC
is defined as
its set of
 proper axioms. 
\item
The second component 
 ``$\, D \,$''
  of a 
 GenAC,  called its 
 {\bf Deductive Apparatus},
is
defined as
the union of its
 logical axioms ``$\,L_D$'' with its
        rules for obtaining inferences.

\end{enumerate}
\begin{example}
\label{nex-2.1}
\rm
This notation 
allows us to
separate  the logical axioms
$~L_D~,~$ associated with  $(  \alpha  , D  )~$, from 
its
``basis axioms'', denoted as ``$\, \alpha \,$''.
It also allows us to  compare
different
%
%
deductive apparatuses
from the literature.
Thus,
the
  $~D_E~$
 apparatus,
from
  Enderton's textbook \cite{End}, 
uses  only  modus ponens
as a rule of inference,
but it deploys a
complicated
4-part  schema of logical axioms.
This differs from
the  $~D_M \,$ and   $\, D_H~$
 apparatuses
in the
Mendelson \cite{Mend} 
 and H\'{a}jek-Pudl\'{a}k \cite{HP91} textbooks.
(They
used  a more reduced set of logical axioms
but employed  ``generalization''
as a second
rule of inference.)
In contrast, the  $~D_F~$ apparatus, from
 Fitting's  and Smullyan's
textbooks \cite{Fi96,Sm95}, 
 uses {\it no logical axioms,}
but  employs a broader
``tableau style'' rule of
inference.
{\bf AN IMPORTANT  POINT} is that while proofs have
different lengths under
different
 apparatuses,
all the common  apparatuses  
 produce  the same  set of  final theorems
from an  initial  common ``axiom basis'' of $\alpha$
(as 
 footnote \footnote{ This is
because all the  common  apparatuses  
satisfy 
the requirements of
G\"{o}del's Completeness Theorem.} explains).
%
%
\end{example}

\begin{definition}
\label{ndef-2.2}
\rm
Let
$ \, \alpha \, $  again 
denote an axiom basis, 
 $ \, D \, $ 
designate
 a
deduction apparatus, and
 $(  \alpha  , D  )$ denote their GenAC.
Henceforth,
the configuration
 $(  \alpha  , D  )$ 
will
be
 called
  {\bf Self-Justifying}
 when
\begin{description}
  \item[  i.   ] one of  $~(  \alpha  , D  )$'s  theorems
(or possibly one of $\alpha$'s axioms)
states that the deduction method $ \, D, \, $ applied to the
basis
system $ \, \alpha, \, $ 
produces a consistent set of theorems, and
\item[  ii.   ]
     the GenAC formalism $ \,( \alpha,D)  \, $ is
actually, in fact,
 consistent.
\end{description}
\end{definition}

\begin{example}
\label{nex-2.3}
\label{ex2}
\rm
Using 
Definition \ref{ndef-2.2}'s 
 notation, our
prior
 research  
\cite{ww93,ww1,ww5,wwapal,ww9}
constructed
GenAC pairs  
$~(  \alpha  , D  )$
that 
 were
``Self Justifying''.
We
also 
proved that
the Incompleteness Theorem 
implies specific
limits beyond which 
self-justifying
formalisms
simply
 cannot transgress.
For any  $\,(\alpha,D) \,$, 
all our articles observed it was
easy
to construct a 
system $ \, \alpha^D \, \supseteq  \,  \alpha  \, $
 that  satisfies
the
Part-i 
condition
(in an isolated context {\it where the Part-ii condition is
 not also
satisfied}).
In essence,
  $ \, \alpha^D \, $  could
consist of all of $~\alpha \,$'s axioms plus 
the added {\bf $\,$``SelfRef$(\alpha,D)$''$\,$} sentence,
defined below: 
\begin{quote} 
$\oplus~~~$ 
There is no proof 
(using 
$D$'s deduction method)
of  $0=1$
from the  {\it union}
 of
the
 axiom system $\, \alpha \, $
with {\it this}
sentence  ``SelfRef$(\alpha,D) \,$'' (looking at itself).
\end{quote}
Kleene 
 \cite{Kl38}
was the first to
notice
how
to
encode
 analogs of 
SelfRef$(\alpha,D) \,$'s  above statement,
which we often 
 call an 
 {\bf $\,$``I AM CONSISTENT'' 
 axiom.}
Each of
Kleene, 
Rogers and Jeroslow 
 \cite{Kl38,Ro67,Je71}
 emphasized
$\alpha ^D$ 
may
be inconsistent
(e.g. 
violate
 Part-ii of   self-justification's
definition
{\it despite}
the assertion in
 SelfRef$(\alpha,D)$'s 
particular
statement).
This is because if the 
 pair $(\alpha,D)$ is too strong
then a
quite conventional
G\"{o}del-style diagonalization argument can
be applied to the axiom basis of
$\alpha^D~~=~~ \alpha \, + \, $ SelfRef$(\alpha,D), ~$
where the added presence of the statement 
SelfRef$(\alpha,D)$ 
will cause this extended version of 
$\, \alpha\,$, ironically,
 to
 become automatically inconsistent.
Thus, an
encoding for
``SelfRef$(\alpha,D)$'' is relatively easy,
via an application of the Fixed Point Theorem,
but this sentence
 is {\it
potentially
devastating.}
\end{example}



\begin{definition}
\label{ndef-2.4}
\rm
Let
 $Add(x,y,z)$ and    $Mult(x,y,z)$ 
denote two 3-way predicates  specifying 
 $x+y=z$ and $x*y=z$.
(Obviously,
   arithmetic's classic
associative, commutative, identity and distributive 
axioms
will 
 have  $\Pi_1$
encodings when
they are expressed  
using these two predicates.)$~~$ 
We will say
that a formalized
axiom basis system of
$~\alpha~$
{\bf recognizes} successor,  addition  and multiplication
as {\bf Total Functions} iff 
it can  prove all of
\eq{totdefxs} - \eq{totdefxm}
as theorems:
\end{definition}
{ \small
\baselineskip =  .9 \normalbaselineskip 
\beq 
\label{totdefxs}
\forall x ~ \exists z ~~~Add(x,1,z)~~
\enq
\beq 
\label{totdefxa}
\forall x ~\forall y~ \exists z ~~~Add(x,y,z)~~
\enq
\beq 
\label{totdefxm}
\forall x ~\forall y ~\exists z ~~~Mult(x,y,z)~
\enq }

\noindent
We will call the GenAC system $(\alpha,D)$ a
{\bf Type-M} 
formalism
iff it proves
\eq{totdefxs} - \eq{totdefxm}
as theorems, {\bf Type-A} if it proves
only \eq{totdefxs} and \eq{totdefxa},
and
it will be called
 {\bf Type-S} if it proves
only \eq{totdefxs} as a
 theorem. 
Also,
 $(\alpha,D)$ 
 will be 
called 
{\bf Type-NS} iff it  can prove
none of \eq{totdefxs} - \eq{totdefxm}. 



\begin{remark}
\label{nremm-2.5}
\rm
The separation of GenAC systems into the
categories of Type-NS, Type-S, Type-A and Type-M systems
helps
 summarize the prior literature
about generalizations and boundary-case exceptions
for the Second Incompleteness Theorem. This is because:
\bed
\item[   $~~~~$i.$~~$  ]
The combined research of Pudl\'{a}k, Solovay, Nelson and Wilkie-Paris
\cite{Ne86,Pu85,So94,WP87},
as  formalized by Theorem $\, ++ \,$,
implies that no
natural Type$-$S  system $(\alpha,D)$
can recognize  its own  consistency
when $D$ represents one of
Example \ref{nex-2.1}'s three
  Hilbert-Frege 
deductive methods
of 
$\, D_E \,$, $\, D_H \,$  
and
$\, D_M$. 
It thus establishes the following result:
\medskip
\begin{quote}
\normalsize \baselineskip = 1.0 \normalbaselineskip 
{\bf ++ }
{\it 
$~~~$
(Solovay's  
modification
\cite{So94}
of Pudl\'{a}k \cite{Pu85}'s formalism 
using some of 
Nelson and Wilkie-Paris \cite{Ne86,WP87}'s
methods)} :
Let 
$ \, (\alpha,D) \, $ 
denote 
a 
Type-S
GenAC system
which assures
the successor operation
will
provably
satisfy
both 
 $  \,   x'     \neq 0     $ and
$     x'     =     y' \Leftrightarrow x=y $.
$~$Then
$ \, (\alpha,D  ) \, $  
cannot verify its own
consistency
whenever
simultaneously
 $D$ is some type of
a Hilbert-Frege
deductive
apparatus and
$~\alpha~$
 treats addition and multiplication
as 3-way relations, 
satisfying 
their usual 
associative, commutative, 
 distributive 
and identity 
axioms.
\end{quote}
\medskip
Essentially, Solovay \cite{So94} 
privately communicated 
to us 
in 1994
an analog of theorem $++$.
Many authors
have noted Solovay
 has 
been
reluctant to publish
his 
nice 
privately communicated
results 
on many occasions
\cite{BI95,HP91,Ne86,PD83,Pu85,WP87}. 
Thus,
approximate  analogs of 
 $++$
 were  explored 
subsequently
 by  Buss-Ignjatovi\'{c},
H\'{a}jek 
and
\v{S}vejdar in \cite{BI95,Ha7,Sv7},
as well as in Appendix A of 
our paper
\cite{ww1} 
and in \cite{wwlogos}.
Also, 
Pudl\'{a}k's initial 1985 article  \cite{Pu85} 
captured
the majority 
of $++$'s 
essence, chronologically before Solovay's observations.
Also,
Friedman did some
closely
related work
 in
\cite{Fr79a}.


\item[   $~~$ii.$~~$  ]
Part of what makes  $++$ interesting is that 
\cite{ww1,ww5,wwapal}
presented  two 
types
of self-justifying
GenAC systems,
 whose
natural hybrid is  precluded by $++$.
Specifically,  these results involve using
Example \ref{nex-2.3}'s 
self-referencing  {\it ``I am consistent''}
 axiom (from 
statement  $\oplus$ ). 
Thus, they establish that
some (not all)
  Type-NS  
systems \cite{ww1,wwapal}
can verify their   own consistency under 
a Hilbert-Frege style deductive apparatus
\footnote{ The Example \ref{nex-2.1} had
provided
three examples of
  Hilbert-Frege style 
deduction operators, called 
$\, D_E \,$, $\, D_H \,$  
and
 $ D_M ~~$. It explained how these 
 deductive operators differ from a tableau-style
deductive apparatus by containing a modus ponens rule.

},
and
some (not all)
 Type-A 
 systems \cite{ww93,ww1,ww5,ww6} can,
likewise, 
corroborate
their 
consistency
under a more restrictive semantic
tableau
 apparatus.
Also, we observed in  \cite{ww2,ww7} how one could
refine $++$ with Adamowicz-Zbierski's
methods
\cite{AZ1} to show 
 most Type-M  systems
cannot recognize their 
semantic tableau
 consistency.
\ennd
\end{remark}

\begin{remark}.
\label{rem2}
\rm
Several of our papers,
starting with our 1993 article \cite{ww93},
 have
used
Example \ref{ex2}'s
 {\it ``I am consistent''} 
axiomatic declaration
$~\oplus~$
 for evading
 the Second Incompleteness Effect.
Other  possible
types of 
            evasions
rest on
     the cut-free methods of
Gentzen and Kreisel-Takeuti \cite{Ge36,KT74},
an
interpretational approach
(such as what
Adamowicz, Bigorajska,
 Friedman, Nelson, Pudl\'{a}k and Visser had
 applied in
\cite{AB1,Fr79b,Ne86,Pu85,Vi5}), or
Artemov's Infinite-Range perspective \cite{Ar19}
(where an infinite schema
  of theorems  replaces
            one 
single
unified consistency
theorem).
We encourage the reader to examine all these articles,
each of which
has 
their own separate 
merits.
Our focus, in this paper,
will be primarily on the next section's
Theorems
 \ref{xxx.2a} and \ref{xxx.2b}.
They
 show
 that 
some
partial
{\it (and not full)} evasions
of the Second Incompleteness Effect
can arise
under a semantic tableau deductive
apparatus.
\end{remark}


\section{ Main Theorems and Related Notation}


\label{sect4}


A function $F $
is called {\bf Non-Growth}
when 
$ F(a_1, \ldots , \, a_j) 
\leq  Maximum(a_1 , \ldots ,  \, a_j)$
holds.
Six  examples of  
non-growth functions are:
\bee
\small 
\item
{\it Integer Subtraction} 
(where $~x-y~$ is defined to equal zero when
 $~x \leq y~),~~$
\item
{\it Integer 
Division}
(where $x \div y$ 
equals
$~x~$ when $y=0$, and
it equals $~\lfloor ~x/y ~\rfloor~$ otherwise),
\item
$~Maximum(x,y),$
\item
$~ Log_{ \, \spadesuit \, }(x)~$ 
 which
is an abbreviation for
$~\lceil~$Log$_2(~x+1~)~\rceil~~$ under 
the conventional
%
notation.
(The  footnote 
\footnote{
The H\'{a}jek-Pudl\'{a}k textbook \cite{HP91} uses the 
notation ``$~\mid  x  \mid~$'' 
 to 
designate
 what we shall call ``$~ Log_{ \, \spadesuit \, }(x)~$''
Thus for $~x \geq 1 \,$, 
 $~ Log_{ \, \spadesuit \, }(x)~$
denotes
the number of symbols that will 
encode the number $~x\,$, when it is written 
in  a binary format. }
explains the 
%
significance
of this
concept.)
\item
$\,~Root(x,y) \, =  \, \lceil  \, x^{1/y} \,  \rceil$ ,  and also
\item $~Count(x,j)$  which  designates the number of
physical
 ``1'' bits
that are stored among $    \,  x$'s rightmost $    \, j    \,  $ bits.
\ene
Our papers
used the term 
{\bf Grounding Function} to refer to
these six 
non-growth operations.
Also,
the term
{\bf U-Grounding Function}
 referred to
a function that corresponds to 
either
one of these six
grounding primitives
or the {\it growth-oriented}
functional operations of
 Addition and
{\it Double$(x)=x+x$}.

Our  language $L^*$,
defined in
 \cite{ww5}, 
  was
built
out of the 
eight U-Grounding function operations
plus the primitives of 
``0'', ``1'',
   ``$ \, = \, $''
and ``$ \, \leq \, $''.
This language
 differs
 from a conventional arithmetic
 by {\it excluding} a
formal
 multiplication function
symbol.
(Instead, it treats
 multiplication as a 3-way relation,
via 
the obvious employment of its
 Division primitive.)
This
notation
 leads 
to a
  surprisingly
 strong 
and tempting
evasion of
the Second Incompleteness Effect.


\begin{deff}
\label{new4def}
\rm
In a context where   $~\, t \, ~$ is  any term
in our language $L^*$,
the 
special
quantifiers
used in
 the wffs
$~ \forall ~ v \leq t~~ \Psi (v)~$ and $~ \exists ~ v \leq t~~ \Psi (v)$
will be called {\bf bounded 
quantifiers}.
Also,
any formula in
our language 
 $L^*$, 
all of whose
quantifiers are 
so
bounded, will
 be called
a $\Delta_0^*$ formula.
The  $~\Pi_n^{* }~$ and  $~\Sigma_n^{* }~$ formulae 
are, thus,
defined by
the
 usual rules, 
{\bf EXCEPT}  they
{\bf DO NOT} contain multiplication function symbols. 
These
rules are that:
\bee
\item
Every  
$\Delta_0^*$ formula will also be a
``$~\Pi_0^{* }~$''  and 
 ``$~\Sigma_0^{* }~$'' formula. 
\item
A
wff
will be  called
 $ \,\Pi_n^{* } \,$
when it is encoded as 
$\forall v_1 ~ ...~ \forall v_k ~ \Phi$  with
$\Phi$ being  $\Sigma_{n-1}^{* }$.
\item
A wff
will be  called
 $\Sigma_n^*$
when it is encoded as 
$\exists v_1 ..\exists v_k ~ \Phi,$  with
$\Phi$  being   $\Pi_{n-1}^{* }$.
\ene
\end{deff}

\begin{remm}
\label{rem-4.2}.
 \rm
A
  sentence $\Psi$ will be called $~$ {\bf  Rank-1* } $~$ when it
can be encoded as either a  $\Pi_1^*$
or  $\Sigma_1^*$ sentence.
 Our definitions for
 $\Pi_1^*$
or  $\Sigma_1^*$ formulae will differ
from Arithmetic's conventional
 counterparts
 by excluding
 multiplication
function symbols.  (This issue will
turn out to 
 be
central
  to our
evasions
of the Second Incompleteness Effect.)
%
\end{remm}

There will be three variants of 
formal deductive
 apparatus
 methods,
which
we will
  compare.
The first is
{\it semantic tableau}.
It will
receive an abbreviated name of
 {\bf ``Tab''}
and correspond to 
Fitting's 
textbook formalism   \cite{Fi96}.
(Its definition can also be found in
the    
      attached
Appendix.)
Thus, a 
Tab-proof for a theorem $~\Psi$,
 from an axiom basis
$~\alpha , \,$ is a
  tree-structure that begins with the sentence
 $~\neg ~ \Psi$ 
 stored inside the tree's root and whose every root$-$to$-$leaf
path establishes a contradiction by containing some pair of contradictory
nodes  that will ``close'' its path. The rules for generating internal
nodes,
 along each
 root$-$to$-$leaf
path,
are that each  node
must be
{\it either} a
proper axiom of $\alpha$ 
$\, or \, $  a deduction
from an ancestor node via one of the  Appendix's 
stated
``elimination''
rules
for the
 $\wedge$, $\vee$, $\rightarrow$, $\neg$,
$\forall$, and $\exists$   symbols.


\smallskip

Our second 
explored
 deductive apparatus
is called {\it Extended Tableau},
and shall be
abbreviated as
{\bf ``Xtab}''. Its definition 
is
 identical to 
 {\bf Tab-}deduction,  except that for any sentence
$\phi$
in our  language $L^*$, the sentence $\phi \, \vee \, \neg \phi$
is allowed 
as an
 internal node in an Xtab proof tree.
(In other words,    {\it Xtab-}deduction
 differs from  {\it Tab-}deduction by allowing all instances
of the Law of Excluded Middle to appear as 
permitted
logical  axioms.
In contrast,  {\it Tab-}deduction 
will  view
 these instances 
   only as
  derived theorems.) 


Our third 
 deductive apparatus
was called
{\bf Tab-1}
in
 \cite{ww5}.
It
is, essentially,
a
  compromise between  Tab and Xtab,
where a 
 ``Tab-1''
 proof for $\Psi$ from
an axiom basis 
$\alpha$ 
corresponds to  a set of
ordered pairs $ \, (p_1,\phi_1), \, (p_2,\phi_2), \, ...., (p_k,\phi_k) \, $
where
\bee
\item $ ~ \phi_k ~  = ~  \Psi \,$ 
\item
Each  $~p_j~$ is a Tab-proof of
what we have called
 a Rank-1* sentence
$~\phi_j~$ from the union of $~\alpha~$ with the preceding Rank-1* sentences
of   $~\phi_1,~\phi_2,~....,~\phi_{j-1}~$.
\ene 
The Rank-1* constraint 
( defined 
 by
Remark \ref{rem-4.2} and utilized by
the above
 Item 2) 
is
%
 significant.
This is because 
 Tab-1 deduction is
 less efficient
 than Xtab 
when {\it the$ \,$former
requires} $\phi_j$ be a Rank-1* sentence. (In contrast, Xtab
{\it  does not impose}
 a similar Rank-1* 
requirement
 upon
$\, \phi \,$ when
 its
Law of the Excluded Middle
 allows  $ \, \phi \, \vee \, \neg \phi$
to appear {\it anywhere}  as a
permissible logical
 axiom,
 {\it for fully arbitrary}  $~\phi.~)$ 
Thus, 
 Xtab
 is
 more desirable than
 Tab-1 {\it when 
it can
actually 
 be
 feasibly (?)
employed.}

Let us say
an axiom system $\alpha$ owns a
{\bf Level-1} appreciation of its own self-consistency
(under a deductive apparatus $D$)
iff it can verify that 
$D$ produces
 no two simultaneous
proofs for a 
 $\Pi_1^*$ 
sentence and 
its negation.
Within this
 context, where $~\aaa~$ denotes any
basis
 axiom
system using
 $L^*\,$'s 
U-Grounding language, IS$_D(\aaa)$ 
was defined 
in \cite{ww5}
to be an axiomatic
formalism  capable of recognizing all of 
$\aaa$'s $\Pi_1^*$ theorems and 
corroborating 
its own Level-1 consistency under $D$'s deductive 
apparatus.
It 
consists
 of the 
following four
groups of axioms:
\begin{description}
\item[Group-Zero:]
Two of the Group-zero axioms will define
the
constant-symbols,
$\bar{c}_0$
and $\bar{c}_1$,
designating the integers of 0 and 1.
The
Group-zero axioms
will
also define the
growth functions of Addition and 
$~Double(x) \, = \, x+x.~$
(They will  enable our formalism to
define
any 
integer
$~n \geq 2~$ 
using fewer than
$3 \cdot \lceil \, $Log$~n~ \rceil \,$ 
logic symbols.)

\item[Group-1:] 
This axiom group  will consist of a
finite set of $\Pi_1^{*} $ sentences, denoted as $~F~$, which
can prove any $\Delta_0^*$ sentence that
holds true under the standard model of the natural numbers.
(Any finite set of 
$\Pi_1^{*} $ sentences $ \, F$, 
with this property,
may be used to define Group-1,
as    \cite{ww5} had  noted.)

\item[Group-2:]
Let $\ulxyz \Phi \urxyz$ denote
$\Phi$'s G\"{o}del Number, and
HilbPrf$_\aaa(\ulxyz \Phi \urxyz,p)$ denote a
$\Delta_0^{*} $ formula indicating that 
$~p~$ is a
Hilbert-Frege styled proof of 
theorem $~\Phi~$ from
axiom system $\aaa$.
For each $\Pi_1^{*}  $ sentence  $\Phi$, the 
Group-2 schema will contain the below axiom 
\eq{group2}.
(Thus IS$_D(\aaa)$ can trivially prove
 all $\aaa$'s 
$\Pi_1^{*}  $ theorems.) 
\begin{equation}
\forall ~p~~~\{~ \mbox{HilbPrf$_\aaa(\ulxyz \Phi \urxyz,p)$}
 ~~
\Rightarrow ~~ \Phi~~\}
\label{group2}
\end{equation}
\item[Group-3:]
The final part of  IS$_D(\aaa)$ 
will 
be a
self-referencing
$\Pi_1^*$
axiom,
that indicates
IS$_D(\aaa)$
is
``Level-1 consistent'' 
under
 $D$'s deductive 
apparatus.
It 
thus amounts to
 the following  declaration:
\begin{quote} 
\# $~${\it No two
proofs exist
for 
a $\Pi_1^{*} $ sentence
and its negation, when
$D$'s deductive 
apparatus
 is applied to an axiom system,
consisting of  
the {\it union}
of 
Groups 0, 1 and 2 with {\bf $\,$this sentence$\,$}
(looking at itself).}
\end{quote}

\smallskip
One encoding for \#
as a self-referencing
$\Pi_1^{*} $
axiom,
had appeared
 in 
\cite{ww5}.
Thus, 
\el{group3} 
is
a $\Pi_1^{*}\, $ 
representation for
  \#
$~$when:
\bed
\item[   a. ] 
$ \mbox{Prf} \, _{\mbox{IS}_D(\aaa)}(a,b) \, $ is
a 
 $\Delta_0^{*} $ formula 
indicating 
that 
$ \, b \, $ is a proof
 of a theorem $\, a\,$
from
the
 axiom basis
  $\mbox{IS}_D(\aaa)$
under
 $D$'s 
deductive apparatus,
$\,~$and
\item[   b. ]  
 $~$Pair$(x,y)$ is a $\Delta_0^{*} $ formula
indicating 
that $ \, x \, $ is  a
 $\Pi_1^{*} $ sentence 
and
 $ \, y \, $ 
represents 
$ \, x \,$'s negation.
\end{description}
\end{description}
\begin{equation}
\forall  ~x~\forall  ~y~\forall  ~p~\forall  ~q~~~~ \neg ~~
[~~ \mbox{Pair}(x,y)~ \wedge ~ 
~\mbox{Prf}~_{\mbox{IS}_D(\aaa)}(x,p)~
\wedge ~ ~\mbox{Prf}~_{\mbox{IS}_D(\aaa)}(y,q)~ ]
\label{group3}
\end{equation}
For the sake of brevity, we will not provide
exact details
about how  
\el{group3} 
can be encoded under the Fixed Point Theorem. 
Adequate 
 details
 are 
provided 
in
 \cite{ww1,ww5}.

\begin{definition}
\label{def3}
Let
 ``$D$'' denote 
any one of 
 the {\it Tab,  Xtab} or
{\it Tab-1}
 deductive
apparatus.
Then we will  say that 
the resulting mapping of  $\mbox{IS}_D(~\bullet~)$ 
 is
{\bf Consistency Preserving}
 iff
 $\mbox{IS}_D(\aaa)$ 
is automatically consistent 
whenever
 all the axioms of $\aaa$ hold
true under the standard model of the natural numbers. 
\end{definition}

The preceding definition raises  questions
about whether the
 mappings of 
    $\mbox{IS}_{\it Tab}(~\bullet~)$, 
  $\mbox{IS}_{\it Tab-1}(~\bullet~)$,
 and  $\mbox{IS}_{\it Xtab}(~\bullet~)$
are consistency preserving.
It turns out that
Theorem \ref{xxx.2a}
will show
the first two of these  mappings 
are  consistency preserving,
while
Theorem \ref{xxx.2b}
explores how
the 
Law of the Excluded Middle
conflicts with 
 $\mbox{IS}_{\it Xtab}(~\bullet~)$'s Group-3 axiom.

\begin{theorem}
 \label{xxx.2a}  The 
  $\mbox{IS}_{\it Tab-1}(~\bullet~)$  and   $\mbox{IS}_{\it Tab}(~\bullet~)$ 
mappings are consistency preserving. (I.e.  the axiom systems 
 $\mbox{IS}_{\it Tab-1}(\aaa )$
and  $\mbox{IS}_{\it Tab}(\aaa )$
are automatically
consistent whenever all $\aaa$'s axioms hold
true
 under the standard model
of the Natural Numbers.)
\end{theorem}

\begin{theorem}
\label{xxx.2b}  In contrast,
  $\mbox{IS}_{\it Xtab}(~\bullet~)$
 fails to be a     
 consistency-preserving mapping.
(More specifically,
  $\mbox{IS}_{\it Xtab}(~\aaa~)$
is
 automatically inconsistent
whenever
$~\beta~$
proves 
some
 conventional $\Pi_1^*$ theorems
stating that
 addition and multiplication
satisfy 
their usual
associative, commutative, 
 distributive 
and identity 
properties.)
\end{theorem}

The proofs of 
Theorems \ref{xxx.2a} and \ref{xxx.2b} 
would be quite lengthy,
if they were derived from first principles.
Fortunately, it is unnecessary for us to do so here
because we gave a detailed justification of 
Theorem \ref{xxx.2a}'s result for  $\mbox{IS}_{\it Tab-1}(~\bullet~)$
in \cite{ww5}, and one can incrementally modify
the 
Remark \ref{nremm-2.5}'s 
special
Invariant of ++ to justify
Theorem \ref{xxx.2b}. Thus, it will be possible for the next two
sections of this paper to 
adequately
summarize
 the intuition behind
Theorems \ref{xxx.2a} and \ref{xxx.2b}, without delving into
  the full
formal
 details.  

Part of the reason Theorems \ref{xxx.2a} and \ref{xxx.2b} are of interest
is because of their
 surprising contrast.
Thus, 
some
 historians have wondered whether
Hilbert and G\"{o}del 
were entirely incorrect 
when their
statements $*$ and $**$
 suggested  some form
of the Consistency Program
would
 likely
be  viable.
Moreover
Gerald Sacks's 
YouTube talk
\cite{YouSa14},
as well as some added comments
by Anil Nerode \cite{Ne20},
have
reinforced this
point.
This is because
 G\"{o}del
repeated analogs of  $**$'s
statement
on several
 occasions,
during 
the later part of his career.
Thus, the
contrast between Theorems \ref{xxx.2a} and \ref{xxx.2b}
provides possible
evidence 
that a
{\it fractional  portion} of what 
 Hilbert and G\"{o}del
had
 advocated,
might become 
feasible.


This
paper
will not have
the page 
 space to go into
the
 full
details, 
but
 the next several sections will
summarize
 the
 gist  behind
the
 proofs
for
 Theorems \ref{xxx.2a} and \ref{xxx.2b}.

\section{Intuition Behind Theorem \ref{xxx.2a}}



Let us recall
the acronym {\bf ``Tab''} stands for semantic tableau
deduction.  This was defined by 
Fitting 
 \cite{Fi72,Fi96}
to be a tree-like proof   
 of a theorem $~\Psi$
 from an axiom basis
$~\alpha~$, whose root consists
of the temporary negated assumption of 
 $~\neg ~ \Psi$ 
and whose every root$-$to$-$leaf
path establishes a contradiction by containing some pair of contradictory
nodes  that ``close''
its 
 path. 
Each internal node along these paths 
must 
{\it either} be a
proper axiom of $\alpha$ 
{\it or be} a deduction
from an ancestor node via one of the 
``elimination''
rules associated with
the
logic symbols of $\wedge$, $\vee$, $\rightarrow$, $\neg$,
$\forall$, or $\exists$ 
(that
are
illustrated in the
Appendix.) 


\begin{example}
\label{ex-5.1}
\rm
Let 
  $\mbox{IS}^M_{\it Tab}(~\bullet~)$ 
denote a mapping transformation identical to
Theorem \ref{xxx.2a}'s
formalism
 of
  $\mbox{IS}_{\it Tab}(~\bullet~)$, 
{\it except that}
  $\mbox{IS}^M_{\it Tab}$
shall  
contain
 a
further  multiplication function operation and,
accordingly,
 have
its Group-3 ``I am consistent'' axiom statements
{\it updated}
 to 
recognize multiplication as a total  function.
It turns out this 
change
will cause
  $\mbox{IS}^M_{\it Tab}(~\bullet~)$ 
to stop satisfying
the consistency-preservation property, which
Theorem \ref{xxx.2a} 
 attributed
 to
  $\mbox{IS}_{\it Tab}(~\bullet~)$. 
\end{example}

\medskip

The intuition behind this 
change
 can be
roughly
summarized
if we let
$   x_0,   x_1,   x_2,     ...    $ 
and  $   y_0,   y_1,   y_2,     ...  $
     denote the
     sequences defined by:
\vspace*{- 0.7 em}
\beq
\label{zs}
x_0~~~=~~~2~~~=~~~y_0
\enq
\beq
x_{i}~~~=~~~x_{i-1}~+~x_{i-1}
\label{as}
\enq
\beq
y_{i}~~~=~~~y_{i-1}~*~y_{i-1}
\label{bs}
\enq
For $\, i \, > \,0 \,$, 
$\,$let $ \, \phi_{i} \, $ 
and $ \, \psi_{i} \, $ 
denote the  
sentences in 
\eq{as} and \eq{bs}
respectively.
Also, 
 let
  $ \, \phi_{0} \, $ and
$ \, \psi_{0} \, $ 
denote \eq{zs}'s
sentence.
Then
 $ \, \phi_0, \, \phi_1, \, ... \, \phi_n \, $
imply
 $ \, x_n \, = \, 2^{ n+1} \, , \, $ and 
 $ \, \psi_0, \, \psi_1, \, ... \, \psi_n \, $
 imply $ \, y_n \, = \, 2^{2^n} \, $.
Thus, the  latter sequence
shall
grow at an
exponentially 
faster
rate than 
the former.
It turns out that this change in growth
speed
causes 
the
  $\mbox{IS}^M_{\it Tab}(~\bullet~)$, 
and
  $\mbox{IS}_{\it Tab}(~\bullet~)$ 
to have 
 quite
 opposite
self-justification
 properties.

In particular,
let
the quantities
Log$(\, y_n \,) \, = \, 2^{n} \, $ 
and
Log$(\, x_n \,) \, = \, {n+1} \, $ 
 represent
the lengths for the binary codings 
for
$ \, y_n \, $
 and 
$ \, x_n \, $.
Thus, $ \, y_n\,$'s 
coding 
will have a length
$\,  2^{n} \,  $, which is
 {\it much larger} than
the $  n+1  $
steps
of
$ \, \psi_0, \, \psi_1, \, ... \, \psi_n \, $
(used
 to
define 
$ \, y_n\,$'s 
 existence).
In contrast,
 $ \, x_n\,$'s 
binary encoding will have a
sharply
smaller length of
size
 $n+1$.
These observations are
significant
because every proof
establishing a variant 
of the Second
Incompleteness Effect
involves 
a 
G\"{o}del
 number $ \, z \,$
encoding 
a  capacity
to self-reference its own definition.

The faster
growing series $y_0,\,y_1,\,,\,...\,y_n$
should,
intuitively, 
have this
self-referencing
capacity because 
 $~y_n\,$'s binary encoding 
has a 
$~2^{n+1}~$  length that 
greatly exceeds
the
size of the $O(n)$
steps 
used
 to define its
value. Leaving aside
many of 
\cite{ww2,ww7}'s
further details,
this
fast growth 
explains
roughly
 why a Type-M 
logic,
such as    $\mbox{IS}^M_{\it Tab}$,
satisfies the semantic tableau version of 
the Second Incompleteness Theorem,
unlike  $\mbox{IS}_{\it Tab}$.


\medskip


Our paradigm also 
explains
why  $\mbox{IS}_{\it Tab}$'s 
 Type-A  formalism
produces
boundary-case exceptions 
for
 the
semantic tableau version of the
 Second
Incompleteness
Theorem. 
This is because  \cite{ww5}
showed that
it was unable
 to
 construct numbers $ \, z \, $ that can
self-reference their own definitions
(when only the {\it  more 
 slowly growing}
              addition
primitive
is available).
In particular
assuming only two bits are needed to encode
each sentence in the sequence
  $   \phi_0,   \phi_1,   ...   \phi_n   $ ,
the length    $     n+1     $ for
 $   x_n   $'s
binary encoding 
is insufficient for encoding this sequence.


Leaving aside many of \cite{ww5}'s details,
this short
length for  $   x_n   $
 explains
the
central
 intuition
 behind
\cite{ww5}'s 
evasion of  
the 
 Second Incompleteness Theorem 
under
  $\mbox{IS}_{\it Tab}~$. 
It arises essentially because of the
{\it sharp}
 difference between the growth
rates of the 
two
 sequences 
of
$~x_1,x_2,x_3 ... ~$ 
and $~y_1,y_2,y_3 ... ~$.


There is obviously insufficient
space for  this extended
abstract
to provide
more
details, here.
A fully detailed
 proof of Theorem \ref{xxx.2a}
 is available
 in \cite{ww5}.
It
 establishes
(see 
\footnote{ 
\baselineskip = 1.4 \normalbaselineskip 
The
{\it exact}
 meaning of this
 implication is
 subtle.
This is 
because Peano Arithmetic (PA)
{\bf  CANNOT$~$KNOW} 
whether
 $\aaa$ is consistent
 when  $ \aaa  \, = \,  PA$.
 Thus,
{\it unlike} 
the 
quite
different
formalism of
 $\mbox{IS}_{\it Tab-1}(PA)~$, the system
of PA   shall linger
in
a state of  self-doubt, about
  whether
both PA and
 $\mbox{IS}_{\it Tab-1}(PA)~$
are
consistent. {\it $~$The {\bf main point} is, however, that we humans
{\bf believe} PA is consistent, and we can use this fact to {\bf confirm} 
that  $\mbox{IS}_{\it Tab-1}(PA)~$ is {\bf BOTH} consistent and able to
verify its self-consistency via its {\bf ``I am consistent''}
axiom. } } )
that
 Peano Arithmetic
can
prove 
 $\aaa$'s  consistency
implies
{\it both }
the consistency and
also the 
 self-justifying property of 
 $\mbox{IS}_{\it Tab-1}(\aaa)$.

Our 
 more 
modest
goal,
within the present abbreviated paper,
has been
to {\it merely }
summarize
the 
 intuition
behind Theorem \ref{xxx.2a}'s 
surprising
 evasion of the Second Incompleteness Effect.       
It arises,
 intuitively, 
because of the
striking
 difference in
the
 growth
rates
between the
two series of
$~x_1,x_2,x_3 ... ~$ 
and $~y_1,y_2,y_3 ... ~$.

\section{Summary of 
Theorem \ref{xxx.2b}'s Proof}

\gvs


A
 formal
 proof of
Theorem \ref{xxx.2b} is
complex,
but it can be 
nicely
summarized.
This is
 because 
 this proposition's
 proof is
similar
to the formal justification for
Remark \ref{nremm-2.5}'s
  Invariant of
 ++ . (The latter's
insight has  come from
 the combined work 
of Pudl\'{a}k, Solovay, Nelson
and Wilkie-Paris
\cite{Pu85,So94,Ne86,WP87}.
It was, also,  subsequently
verified
by   several
 other
authors  \cite{BI95,Fr79a,Ha7,Sv7,ww1}
in  slightly different forms.)

The
crucial aspect of the  Hilbert-Frege
deductive methodology
 is that 
its
modus ponens rule  
 assures thar
a proof of a theorem 
$\psi$ 
from an axiom system $\alpha$ 
has a length
  no more than 
proportional to
the sum of the proof-lengths used to derive
 $\phi$ and $~\phi \rightarrow \psi~$.
 This 
{\bf ``Linear-Sum Effect''
}
 {\it does
not}
 apply
also
to
{\it Tab-}deduction
(because 
the latter
 lacks a modus ponens rule).

The {\it Xtab}
deductive 
methodology is , however, quite
different
from the
{\it Tab} form of deduction,
in that {\it  only Xtab}
supports an analog of the 
prior paragraph's
 ``Linear-Sum Effect''.
This is because
 any node of an Xtab proof-tree
is allowed
to store any sentence of the form $~ \phi \, \vee \, \neg ~\phi~$
(as a consequence of its allowed use of the Law of Excluded
Middle).
This added feature
will 
allow an
 {\it Xtab} proof for
$\psi$ to have a length
proportional to
 the sum of
the proof lengths for 
 $\phi$ and $~\phi \rightarrow \psi~$.
In particular,
such an
 Xtab 
proof
for $~\psi~$
 will
 consist
 of
 the following 
four steps:
\bee
\item
The root of 
an Xtab
 proof for $\psi$
consists of
the usual temporary negated hypothesis of
 $~\neg \psi~$ (which the remainder of the proof tree will
show is impossible to hold). 
\item
The child of this 
root
node 
consists of
an allowed invocation of the
Law of the Excluded Middle of the 
{\it particular}
 form $~ \phi \, \vee  ~\neg \, \phi~$.
\item
The relevant Xtab
proof tree will next employ
the Appendix's branching rule for allowing the two
sibling nodes of
 $~  \phi ~$
 and $~\neg \, \phi~$ to descend from
Item 2's node.
\item 
Finally, our Xtab proof will insert below
(3)'s  left sibling  node
of  $~  \phi ~$ a subtree 
that
 is no longer than
a proof for
   $~\phi \rightarrow \psi~$,
and likewise
insert
 a   proof for $~\phi~$ 
below (3)'s right sibling
of  $~ \neg ~ \phi~$.
\ene


The point is that the
very
 last step of the above 4-part
 proof
has a length no greater than the sum of the two
 proof lengths for 
 $\phi$ and $~\phi \rightarrow \psi$.
(This is analogous to the proof expansions resulting from a
conventional
modus ponens operation.)
Its first three steps
will
 have
{\it entirely 
 inconsequential} effects
that
increase 
the
 overall proof length by no
more than
a {\it   tiny}
amount, that is  
 proportional 
to the
trivial sum
of the lengths for the two individual
sentences of ``$~\phi ~$'' and ``$ ~\psi~$''.


Hence,
the preceeding 
 ``Linear-Sum Effect'' allows us
to construct  an analog 
of Remark \ref{nremm-2.5} 's
earlier  Theorem $\, ++ ~$ 
for Xtab deduction.
It is formalized by the statement $~\bigodot~$
 below:

\bigskip

\begin{description}
\small
 \baselineskip = 1.15 \normalbaselineskip 
\item[ $\bigodot~~$]
Any axiom system $~\cal{A}~$
 is
 {\it  automatically inconsistent} whenever
it satisfies the following three conditions:
\bed
\item[   I. ]
$~~~\cal{A}$
 can verify Successor is a total function (as
 \el{totdefxs}  formalized).   
\item[   II. ]
$~\cal{A}$
can prove
 addition and multiplication
(viewed as 3-way relations)
satisfy 
their usual 
associative, commutative, 
 distributive 
and identity-operator  properties.
\item[ III.]
$~\cal{A} \,$
proves 
an added
 theorem (which turns out to be false) affirming its own
consistency when the
Xtab
 deductive apparatus is used.
\ennd
\end{description}

\bigskip


It is not possible to provide a short  proof for
statement $~ \bigodot~~$
because it will rest upon the very
detailed
 ``Definable Cut'' machinery
from
 pages 172-174 of the
 H\'{a}jek-Pudl\'{a}k textbook \cite{HP91}.
The intuition behind  $~ \bigodot~~$ is,
however,
quite simple.
 It is that
statement  $~ \bigodot~~$ causes
 ++'s 
mechanism
to
generalize  from Hilbert-Frege  
deduction to Xtab (because both satisfy the 
Linear-Sum Effect).

The nice aspect of
 $~ \bigodot~~$
is that its machinery
 establishes
 Theorem \ref{xxx.2b}. 
This is because if $~\beta~$ satisfies
 Theorem \ref{xxx.2b}'s hypothesis
then 
 $\mbox{IS}_{Xtab}(\aaa)$ 
will satisfy\footnote{Actually,  $\mbox{IS}_{Xtab}(\aaa)$
will satisfy a requirement stronger than Item I
because it recognizes addition as a 
total function.}
 the      
conditions I-III that cause  
  $\mbox{IS}_{Xtab}(\aaa)$ 
to become inconsistent.



\section{More Elaborate Forms of Theorems \ref{xxx.2a} and \ref{xxx.2b}}


Our
results 
in  Theorems
 \ref{xxx.2a} and 
\ref{xxx.2b}
demonstrate 
 self-justifying methodologies
apply to
 ``Tab''',
{\it  but
  not} also   ``Xtab''  deduction.
(This is because Xtab
treats the
the Law of Excluded Middle
as  a  formal
schema of logical axioms, and the latter
activates the power 
of the Second Incompleteness Effect.)

Our
 goal
in  this 
section
will be
 to view
this machinery
in more meticulous detail.
 Thus, we 
 will
 explore 
{\it at what
exact juncture}
the boundary
is
 crossed between
 generalizations of
the Second Incompleteness Theorem and
its
permissible
exceptions.

\begin{deff}
\label{def-n1}
\rm
Let
 $L^*$ 
 again denote the base arithmetic language
(that was defined in \textsection \ref{s3} ), and 
 $~Z~$
denote an arbitrary set of sentences appearing in the  language $~L^*$
(such as its set of $\Pi_2^*$ sentences).
Let us 
recall that
the Appendix defined
 a semantic tableau proof of a theorem $\Psi$ from
 $\alpha$'s axiom system.
Then
 a  {\bf Z-Enriched}
 modification for a 
semantic tableau
proof of 
a theorem
$\Psi$, from
$\alpha$'s set of proper axioms,
 will be
 defined 
 as
the particular
 refinement of the Appendix's
proof-tree
 formalism
that 
allows
\el{nneq} as 
an
 added  permissible
 logical axiom,
for
 any
 $ ~\Upsilon~ \in Z$.
\beq
\label{nneq}
\Upsilon  ~~\vee ~~ \neg ~\Upsilon
\enq
\end{deff}

\begin{deff}
\label{def-n2}
\rm
It is also
of
 interest to consider a 
slight modification of the preceeding nomenclature, 
where $Z$ is a set of formulae that are 
allowed to be
free
in the
single
 variable of  $~x~$
(instead of
representing
a  sentence that contains no free variables).
In this case,  $\Upsilon(x)$ will designate a formula, 
within the subset of 
 $Z$, and
\el{nneq2} 
will replace
\el{nneq} as 
the
added permissible
logical axiom
that
can
be allowed to
appear inside a 
{\bf ``Z-Base Variable Enriched''} proof.
\beq
\label{nneq2}
\forall ~x~~~~
\Upsilon(x)~\vee ~~ \neg ~\Upsilon(x)
\enq
\end{deff}

A fully detailed justification will not be 
provided
here,
but it turns out
our results from \cite{ww1}
can be expanded to show that
their
 evasions of 
the semantic tableau version of the Second Incompleteness
Theorem can be extended to both the cases of
Z-Enriched and ``Z-Base Variable Enriched'' 
mechanisms,
when
 $~$Z$~$
 represents 
the
 $\Delta^*_0$ class of formulae. 
We can also extend
our results from \cite{wwlogos}
to show that the comparable
evasions of the semantic tableau version
of the Second Incompleteness Effect 
will fail
at
and above
 the $\Pi_2^*$ level.  

We
 conjecture 
the preceeding 
$\Delta_0^*$
evasions of the Second Incompleteness
Theorem 
will
continue
at the  $\Pi_1^*$ level, but this fact has not 
yet been formally
proven.


A fascinating aspect about this subject
is that
semantic tableau deduction satisfies
its particular variant of 
G\"{o}del's Completeness Theorem
 \cite{Fi96,Sm95}.
 Thus,  the set
of theorems proven
by
an axiom system $~\alpha~$, via
a conventional (unenriched) version of semantic tableau deduction,
{\it is identical to}  those theorems 
proven
by a Z-enriched deductive
mechanism.
Yet
despite this  invariance, the  proof-lengths
change,
{\it quite 
sharply,}
under the  Z-enriched formalisms of 
Defintions 
\ref{def-n1} and \ref{def-n2}.
This 
extreme
change in proof-length
 causes
the deployment of an  
 {\it ``I am consistent''}
axiom
to 
 become
{\it
fully infeasible} 
when $\Upsilon$  in \el{nneq} 
is allowed to represent any arbitrary   $\Pi_2^*$ sentence
(see footnote \footnote{The point is that the sharp compression
in proof lengths  produces G\"{o}del-like Diagonalization compressions,
similar to those
particular
 Second Incompleteness Effects
applicable to
       $\Pi_2^*$  sentences,
that are
examined in \cite{wwlogos}.} ).
\gvs

\section{Further Generalizations}
\gvs

\normalsize \baselineskip = 1.37 \normalbaselineskip

\label{sect8}


For the sake of simplicity, the previous sections
had focused on the semantic tableau  deductive apparatus.
However, 
it is known 
 \cite{Fi96}
that resolution shares numerous characteristics
with tableau.  Therefore, it turns out that 
 Theorems \ref{xxx.2a} and \ref{xxx.2b} do generalize when resolution
replaces semantic tableau.

In particular,
 let us say a theorem $~\cal{T}$ has a
$Res-$proof
from  $~\alpha~\,$'s set of proper axioms
when there is a resolution-based proof
  \cite{Fi96}
 of $~T~$ from  $~\alpha~\,$.
Also,  the term 
$Xres-$proof of $~T~$
 refers
to the obvious extension
of a  $Res-$proof that allows
all instances of
 the Law of Excluded Middle 
(from the base language of $L^*~)$
to appear as formalized logical axioms.

It turns out
$Xres$ 
differs from $Res$ in the same manner
$Xtab$ differed from $Tab$.    
Thus,
the
 obvious generalizations of  Theorems \ref{xxx.2a} and \ref{xxx.2b}
hold for $Res$ and
$Xres$.  
In particular,
   $\mbox{IS}_{\it Res}(~\bullet~)$ 
is a  consistency preserving
transformation, but
   $\mbox{IS}_{\it Xres}(~\bullet~)$ 
again 
is not.



\smallskip

Some logicians may, also, wish to examine
 special
 speculations
 in \cite{ww16}'s  arXiv
article.
It contemplated    an 
 alternative approach, where 
self-justifying 
arithmetics 
  employ 
an 
unconventional
``indeterminate''
 functional
object,
called the
 $\Theta$  primitive,
 to formalize the traditional
 properties of
an 
endless sequence of integers.

If a conjecture
stated in  \cite{ww16}
 is
correct
(as we are almost certain it is),
 then
 such
a self-justifying
machine
will be plausible 
 for constructing
 the entire
 set of natural numbers, 
{\it without encountering} the
usual incompleteness difficulties 
that
the Theorem ++
(of
 Pudl\'{a}k and Solovay)
associated with
Type-S formalisms
(that recognize merely 
Successor as 
a total function).
Interestingly, the
$~\theta~$ function primitive
of  \cite{ww16}
should
 allow a substantial 
 Type-NS arithmetic
 to  exist
that can
 {\it simultaneously}
 recognize its own
Hilbert-Frege consistency 
and
      possess
a formalized ability 
to constructively
enumerate
 the 
full
infinite collection of integers $~0,~1,~2,~3,~ .... ~$


\section{ Ironic Events and
 Related Speculations:}



\gvs

The
initial 19-page draft 
of this article
was accepted by the LFCS-2020 conference
and was published by Springer
 \cite{ww20}, 
 shortly
before
the 
 Covid crisis 
 commenced.
During  January 4-7,
when LFCS met,
there was little knowledge 
\footnote{This is because the Chinese
authorities 
 announced the presence of
Covid only on
 December 31.  Their announcement had not yet
attracted any attention
at the LFCS-2020 conference.}
about  
        the
soon-to-appear 
 epidemic.
The nature of
the Covid
event
did   become apparent
by  March  of 2020.
At that time, the 
ASL
 changed
 its previously  planned North American
Annual Meeting into
a virtual conference (with a
virtual
presentation
of our
planned slides being posted
 at the ASL's
web site).

This ironic chronology
is, perhaps,
worth briefly  
recollecting
because
of the
connection between
\thx{xxx.2a}
  with
the new world of computing that is
now, 
 currently,
emerging.

Thus, mankind will likely become increasingly dependent upon
computers in the future.  
For instance, the spread of a serious epidemic can be 
more effectively
contained, when
  staffs at  medical facilities are computerized, as
much as 
 feasible.
            (Then a patient, suffering
from  a virus, will be less
contagious, 
when
 virus particles bounce off
sterile
 computerized robots,
instead of
            encountering vulnerable  human staff members.)
Also,  transportation
networks and factories processing food-materials
will
be
safer if they are run by
computerized robots,
rather than depend upon human beings
(who
breathe 
out
air containing
dangerous virus
 contagions).


Our point is that a large variety of forms of Artificial Intelligence
will likely 
become increasingly prominent in the future. 
Thus,
 AI-based machines
should become 
more
effective, if their actions are both consistent and 
display the maximal amount of awareness about
their self-consistency
(that is plausible under
anticipated future
 generalizations of
our
 self-justifying
formalism).

It is clear that
 AI-based
computers
will
exhibit  a broad  variety of
automated
skills, many of which
will be
 {\it  only partially related to} 
\thx{xxx.2a}'s self-justifying 
 IS$_{Tab}$ 
mechanism.
(For instance, future 
AI-machines
will,
certainly, 
need automated skills 
that
 master the  arts of
visual
learning, motion planning and
several
forms of
decision-making.)
Nevertheless,
a 
quite
fascinating point is  that
the early 20-th century predictions of
 Hilbert and G\"{o}del, 
in
 $*$ and $**$,
 will gain 
some
new positive interpretations, 
when they 
anticipated
significant
benefits
from future generations of
thinking  machines
  being aware
about the
 consistency of,   at least,
 {\it  their
 specialized restricted} forms of 
mathematical
knowledge.


We do not wish to pursue
these
points
further, here,
 because there will, certainly,  be many other
types  of 
 unanticipated
events, which also 
 advance the need for more elaborate forms
 of Artificial Intelligence in the future.
These
 future events should
be consistent with,
at least,
 the broad 
 predictions that
 Hilbert and G\"{o}del 
made in their
famous
 statements
 $*$ and $**$.

\section{Concluding Remarks}

Our main results  in this article
 are
surprising
   because it is
quite
unusual 
 for an initially consistent formalism
$~\alpha~$
{\it to
become inconsistent} when its 
initial
schema of theorems (establishing
the 
widespread
validity of the Law of the Excluded Middle)  
is
transformed
 into
being a
schema of logical axioms.

This
unusual
 effect 
arose because
the meaning of a
Group-3  {\it ``I am consistent''}
axiom
changes,
{\it 
 quite
substantially,}
when
theorems
 are 
transformed
 into logical
axioms
(as illustrated by footnote 
\footnote{The 
point is that proofs are compressed when theorems
are transformed into  logical axioms, and such  compressions
can produce diagonalizing contradictions under some Type-A
logics using  {\it ``I am consistent''}  
axioms.} ).
Thus,
unacceptable
diagonalizing contradictions
can 
occur
when an 
 {\it ``I am consistent''}  
axiom is able to reference itself
{\it in the context of a
 SUFFICIENTLY
 POWERFUL 
mathematical machine.}

The
 contrast between Theorems \ref{xxx.2a} and \ref{xxx.2b} 
(where only the former 
 eschews
diagonalization effects)
helpfully explains
     how Hilbert and
G\"{o}del
  appreciated the
Second Incompleteness Effect,
{\it while 
they were
 simultaneously
 cautious about it.}
%
%
Moreover,  G\"{o}del's
particular
 remark $**$
 should not be ignored when
comments from
Gerald  Sacks
and 
Stanley  Tennenbaum  \cite{Ne20,YouSa14}
recalled how 
G\"{o}del
reiterated the gist of his
1931-published remark,
many years after its printing.
 Indeed,
 it   
 is      noteworthy
Harvey Friedman
recorded a YouTube 
lecture 
 \cite{Fr14},
stating
            he was
also tentatively
 open to the
possibility that the
 Second Incompleteness
Theorems might allow
 partial exceptions.


Thus, while there is no
doubt that the
Second Incompleteness Theorem
will
 be remembered for its
seminal impact,
$~$ its
part-way exceptions
are also 
significant.
This is because futuristic high-tech computers
will better understand their self-capacities,
if they
own
some {\it  partial}
awareness about their own consistency.




There is no page space to
delve
into all
details
here.
However,
the distinction between 
the initial ``IS(A)'' system, 
from
 our 1993 and 2001 papers
\cite{ww93,ww1}, with the more
sophisticated
 $\mbox{IS}_{\it Tab-1}(\aaa)$
formalism 
of our year-2005 article
\cite{ww5} 
should,
 also, 
 be
 briefly
 mentioned.
Our older ``IS(A)'' formalism was 
actually
simpler, 
 but it was
substantially
weaker because it
 only recognized the
non-existence of a proof of $0=1$ from itself.
In contrast,   $\mbox{IS}_{\it Tab-1}(\aaa)$'s 
 Group-3 axiom 
can corroborate that
{\it no two  simultaneous proofs} exist for a
Rank-1* sentence
and its negation.
This is an important distinction,
because
the First Incompleteness
Theorem
 indicates
no decision procedure 
exists
 for
separating 
all  
true from false Rank-1* sentences.
(See 
\cite{wwapal,ww6,ww9,ww16} for other 
particular refinements
for
 our ``IS(A)'' formalism.)



In summary,
the
 main
  purpose of
this article
has been
 to explore the contrast between the opposing
Theorems \ref{xxx.2a} and \ref{xxx.2b}.
The latter  theorem,  thus, 
provides {\it another 
helpful
reminder }
about
the millennial importance of 
G\"{o}del's seminal
 Second
 Incompleteness Theorem.
Yet at the same time,
Theorem \ref{xxx.2a} 
illustrates how
some {\it partial  exception}s to G\"{o}del's result
do arise, as Hilbert and G\"{o}del
predicted   in their
statements $*$ and $**$.





In essence,
the 2-way contrast between Theorems \ref{xxx.2a} and \ref{xxx.2b}
may
 be as 
significant
 as their
 individual
actual
 results.
%
%
This is because
 the Second Incompleteness Theorem is
fundamental
 to
 Logic.
 Many historians have,
 thus,
been
perplexed by
the
{\it partial} reluctance that 
 Hilbert and G\"{o}del
had expressed about it in
 $*$ and $**$. 
A  
partial
reason for this 
reluctance
is, perhaps, 
related
 to the
contrast
 between 
these two opposing theorems.


\smallskip

{\bf ACKNOWLEDGMENTS:} I thank Seth Chaiken 
and James P. Torre, IV
for several quite helpful comments about
how to
improve the presentation.

  \newpage

\section*{Appendix providing a
formal
 definition for a Semantic Tableau proof:}

\gvs

Our definition of a semantic tableau proof 
is similar to 
analogs
from  the textbooks by
Fitting and Smullyan 
 \cite{Fi96,Sm95}.
A tableau proof  of a theorem $\Psi$ from a set of proper
axioms 
(denoted as $~\alpha~$)
is therefore
 a tree structure, whose 
root contains the temporary contradictory assumption of $~\neg \, \Psi~$
and whose every descending root-to-leaf branch affirms a contradiction
by containing both some sentence $\phi$ and its negation  $\neg \, \phi$.
Each internal node in this tree will be either a proper axiom of
$~\alpha~$ or a deduction
from
a higher 
ancestor
 in this tree via one of
six 
elimination
rules for the logical connective symbols of
$~\wedge~$,  $~\vee~$,   
$~ \rightarrow ~$, $~\neg~$, $~\forall~$ and $~\exists~$.
(These  rules
use a notation where
``{\bf $ \,  $A$  ~  \Longrightarrow  ~  $B$  \, $}''
is an abbreviation for a sentence
 {\bf $ \,   $B$  \, $} 
being
 an allowed deduction
from its ancestor of  {\bf $ \,  $A$  ~~ $}.)
\begin{enumerate}
\itemsep 5pt
\normalsize
\item $~ \Upsilon \wedge \Gamma \, ~ \Longrightarrow ~ \, \Upsilon
~$ 
and 
$~ \Upsilon \wedge \Gamma \, ~ \Longrightarrow ~ \, \Gamma ~$ .
$~~~$
\item $ \,  \neg  \,\neg \, \Upsilon  \,  \Longrightarrow  \,  \Upsilon. \, $  
Other rules for
the ``$ \, \neg \,$'' symbol  are:
$ \, \neg ( \Upsilon \vee \Gamma )  \,  \Longrightarrow  \,  \neg \Upsilon
\wedge \neg \Gamma$,
$ \, \neg ( \Upsilon \rightarrow \Gamma )  \,  \Longrightarrow  \,   \Upsilon
\wedge \neg \Gamma \, $,
$ ~~~~\, \neg ( \Upsilon \wedge \Gamma )  \,  \Longrightarrow  \,  \neg
\Upsilon \vee \neg \Gamma \, $,
 $~ \,   \neg \, \exists v \, \Upsilon  (v)  \,  \Longrightarrow  \,  
\forall v   \neg \, \Upsilon  (v)  \, $
 and $ ~~\,   \neg \, \forall v \, \Upsilon  (v)  \,  \Longrightarrow  \,  
\exists v \,  \neg  \Upsilon  (v)$
\item A pair of sibling nodes $~ \Upsilon ~$ and $~ \Gamma ~$ is
allowed
when their ancestor is
$~\Upsilon \, \vee \, \Gamma.~$
\item A pair of sibling nodes $ \,  \neg \Upsilon  \, $ and $ \,  \Gamma  \, $ is
allowed
when their ancestor is
$ \, \Upsilon \, \rightarrow \, \Gamma$.
\item $\forall v \, \Upsilon  (v)  \,  \Longrightarrow    \, \Upsilon(t)  \, $
where $t$  may denote any 
term.
\item $~ \exists v \, \Upsilon  (v) ~ \Longrightarrow ~ \, \Upsilon(p) ~$
where $\,p \,$ is a newly introduced parameter symbol. 
\end{enumerate}
One minor difference in  notation is  we treat 
 ``$~ \forall~ v \leq s~~~ \Phi(v)~$''
as an abbreviation for
 $~ \forall  v  ~~ \{ ~ v \leq s~~ \rightarrow ~~\Phi(v)~ \}~$
and
 ``$~ \exists~ v \leq s~~~ \Phi(v)~$''
as an abbreviation for
 $~ \exists  v  ~~ \{ ~ v \leq s~~ \wedge ~~\Phi(v)~ \}~$.
Therefore,
Rules  5 and 6 
imply the following
 hybrid 
rules
for processing
 bounded universal and
 bounded
 existential quantifiers:
\begin{description}
\itemsep 5pt
\normalsize
\item[  a. ]  
 $\forall v \leq s  \, \Upsilon  (v) ~~ \Longrightarrow ~~
t \leq s \, \rightarrow \, \Upsilon(t) $ 
where $\,t \,$ may be any arithmetic term.
\item[  b. ] 
 $~ \exists v \leq s ~ \, \Upsilon  (v) ~~ \Longrightarrow ~ ~
p \leq s ~ \wedge~ \Upsilon(p) ~$
where $\,p \,$ is a new parameter symbol. 
\end{description}

{\bf Added Comment:}
The preceding paragraph has formalized what \textsection \ref{sect4}
called the ``Tab'' version of a semantic tableau proof.
Its ``Xtab'' variant is identical
 except that  any node
may optionally 
store
a sentence of the form 
$ \, \mho  ~\vee ~ \neg \, \mho~$ 
(for arbitrary $\mho~)$,
as a manifestation of 
its
%
 allowed use of the 
Law of the Excluded Middle.


\end{document}